\documentstyle[amssymb,11pt]{article}
\title{\bf Derivation algebras of toric varieties}
\author{Antonio Campillo, Janusz Grabowski, and Gerd M\"uller}
\date{January 1997}
\begin{document}
\newcommand{\sub}{\subseteq}
\newcommand{\F}{{\cal F}}
\newcommand{\Der}{{\rm Der}\:}
\newcommand{\rk}{{\rm rk}\:}
\newcommand{\D}{{\Bbb D}}
\newcommand{\St}{\tilde{S}}
\newcommand{\T}{\Theta}
\newcommand{\TS}{\T(S)}
\renewcommand{\l}{\lambda}
\newcommand{\Tl}{\T_\l}
\newcommand{\tl}{t^\l}
\newcommand{\m}{\mu}
\newcommand{\Tm}{\T_\m}
\newcommand{\tm}{t^\m}
\newcommand{\Tlm}{\T_{\l+\m}}
\newcommand{\Dl}{\D_\l}
\newcommand{\Dli}{D_{\l i}}
\newcommand{\Dlj}{D_{\l j}}
\newcommand{\Dmj}{D_{\m j}}
\newcommand{\Yli}{Y_{\l i}}
\newcommand{\n}{\nu}
\newcommand{\mi}{\m_i}
\newcommand{\lj}{\l_j}
\newcommand{\li}{\l_i}
\newcommand{\lm}{\l_m}
\newcommand{\Dlmi}{D_{\l+\m,i}}
\newcommand{\Dlmj}{D_{\l+\m,j}}
\renewcommand{\L}{\Lambda}
\renewcommand{\d}{\partial}
\newcommand{\di}{\partial_i}
\newcommand{\dj}{\partial_j}
\newcommand{\Z}{{\Bbb Z}}
\newcommand{\N}{{\Bbb N}}
\renewcommand{\P}{{\Bbb P}}
\newcommand{\A}{{\Bbb A}}
\newcommand{\Zn}{\Z^n}
\newcommand{\Nn}{\N^n}
\newcommand{\G}{\Gamma}
\newcommand{\DE}{\Delta}
\newcommand{\ad}{{\rm ad}\:}
\renewcommand{\a}{\alpha}
\newcommand{\ai}{\a^i}
\renewcommand{\b}{\beta}

\maketitle
\section{Introduction}
Normal affine algebraic varieties in characteristic 0
are uniquely determined (up to
isomorphism) by the Lie algebra of derivations of their coordinate
ring. This was shown by Siebert \cite{Si} and, independently, by Hauser
and the third author \cite{HM}. In both papers the assumption of
normality is essential. There are non-isomorphic non-normal varieties
with isomorphic Lie algebras.
The third author \cite{M} treated certain non-normal varieties
defined in combinatorial terms by showing that closed simplicial
complexes can be reconstructed from the Lie algebra of their
Stanley-Reisner ring.
Here we study this problem for
(in general, non-normal) toric
varieties defined by simplicial affine semigroups.
\\[1ex]
We show that such
toric varieties are uniquely determined by their Lie algebra if they
are supposed to be Cohen-Macaulay of dimension $\ge2$. The
corresponding statement is false in dimension 1. For toric
curves we need
the stronger hypothesis that they are Gorenstein.
In fact, we can reconstruct from the Lie algebra
the semigroup defining the variety.
Our result should be compared with a recent one of
Gubeladze \cite{Gu} saying that an affine semigroup is uniquely
determined by the toric variety it defines (more precisely, by its
coordinate ring as an augmented algebra).
\\[1ex]
The main tool in our proofs is a root space decomposition of the Lie
algebra of derivations of a Buchsbaum semigroup ring. The set of
roots is closely related to the underlying semigroup. This structural
description will be used to prove two more results. We show, in the
Cohen-Macaulay case, that every automorphism of the Lie algebra is
induced from a unique automorphism of the variety. And we establish
an infinitesimal analogue of the last statement: Every derivation of
the Lie algebra is inner, i.~e., the first cohomology of the Lie algebra
with coefficients in the adjoint representation vanishes.
\\[1ex]
Our results were obtained during visits at the
Mathematics Departments
of the Universities in Valladolid, Warszawa, and Mainz.
We thank these institutions (as well as the Spanish-German
Acciones Integradas) for financial support and their members for
their hospitality.

\section{The root space decomposition}
Let $S$ be an affine semigroup, i.~e., a
finitely generated subsemigroup of some $\Nn$.
We stress that, in this paper, semigroup always means semigroup with zero
element.
Denote by $G=G(S)$
the subgroup of $\Zn$ generated by $S$ and by $r=\rk S=\rk G(S)$ its rank.
Let $C_S$ be the convex polyhedral cone spanned by $S$ in ${\Bbb Q}^n$.
We shall suppose throughout that $S$ is {\em simplicial}, i.~e., that the
convex cone $C_S$ can be spanned by $r$ elements of $S$.
For an algebraically closed field $k$ of
characteristic 0 let $k[S]\sub k[t]=k[t_1,\ldots,t_n]$ denote the
corresponding semigroup ring. We need to recall how the property of
$k[S]$ being Cohen-Macaulay or Buchsbaum can be described in terms of $S$.
For this purpose, let $F_1,\ldots,F_m$ be the
$(r-1)$-dimensional faces of $C_S$.
Set
\[S'_i=\{\l\in G,\;
  \l+s\in S\;\mbox{for some}\;s\in S\cap F_i\}
\]
for $i=1,\ldots,m$, and
$S'=\bigcap S'_i$.
\begin{trivlist}
  \item[] \bf Proposition 1. \it
  For a simplicial affine
  semigroup $S$ the semigroup ring $k[S]$ is Cohen-Macaulay
  (resp.\ Buchsbaum) if and only if $S'=S$
  (resp.\ $S'+(S\setminus\{0\})\sub S$).
\end{trivlist}
For the proof see
\cite{GSW}, \cite[Theorem 6.4]{St},
\cite[section 4]{TH}, and \cite[section 6]{SS}.
The semigroup $S'$ is called the {\em Cohen-Macaulayfication} of $S$.
Let
\[\bar{S}=\{s\in G,\; ms\in S \;\mbox{for some}\; m\in\N,m\not= 0\}.
\]
It is known \cite[section 1]{Ho} that $k[\bar{S}]$ is the normalization
of $k[S]$. An affine semigroup $S$ is called {\em standard} if
  \begin{itemize}
    \item[(i)] $\bar{S}=G(S)\cap\Nn$.
    \item[(ii)] For all $i$ the image of $S$ under the
    the projection $\pi_i$ on the $i$-th component is a numerical
    semigroup, i.~e.,
    the complement $\N\setminus\pi_i(S)$ is finite.
    \item[(iii)] The semigroups $S\cap \ker\:\pi_i$, $i=1,\ldots,n$,
    are distinct of rank equal to $\rk S-1$.
  \end{itemize}
It was shown by Hochster \cite[section 2]{Ho} that every
affine semigroup is isomorphic to a standard one. Hence
we shall assume throughout that $S$ is standard. In that case
the cone $C_S$ has exactly $n$ faces of dimension $r-1$, namely
the convex cones spanned by the $S\cap\ker\:\pi_i$. Hence
\[S'_i=\{\l\in\Nn,\;\l+s\in S\;\mbox{for some}\;s\in S\;\mbox{with}\;
       s_i=0\}
\]
for $i=1,\ldots,n$.
A standard affine semigroup $S$ is simplicial if and only if
$S$ has elements
on every coordinate axis. In fact, the cone of
a simplicial affine semigroup of rank $r$
has only $r$ faces of dimension $r-1$. Standardness gives $r=n$.
Then the edges of $C_S$ are the intersections of $C_S$ with the
coordinate axes, see \cite[section 1]{SS}.
The reversed implication is obvious.
Let $a_i\in\N$, $a_i\not=0$, be the minimal number
such that $\ai=(0,\ldots,0,a_i,0,\ldots,0)\in S$,
where the nonzero entry is at the $i$-th place.
\begin{trivlist}
  \item[] {\bf Proposition 2.}
  \it Every $k$-linear derivation $D$ of $k[S]$ extends uniquely to a
  derivation of the polynomial ring $k[t]$.
  \item[] Proof. \rm As $S\sub\Nn$ is standard
  and simplicial it has rank $n$
  and $k[S]$ has dimension $n$. Hence the rational function field $k(t)$
  is a separable finite extension of the quotient field $k(S)$ of $k[S]$.
  Therefore $D$ extends uniquely to a derivation $D$ of $k(t)$.
  Write $D=\sum f_i \di$ with
  $f_i\in k(t)$, say $f_i=g_i/h_i$ with coprime $g_i,h_i\in k[t]$.
  With the semigroup elements $\ai$ introduced above we have
  \[a_it_i^{a_i-1}f_i=D(t^{\ai})\in k[S] \sub k[t]
  \]
  and $h_i$ divides $t_i^{a_i-1}$. As $\pi_i(S)$ is a numerical semigroup
  there is $s\in G$ with the $i$-th component $s_i=1$. Using simpliciality
  we may assume that $s\in\Nn$, hence $s\in\bar{S}$.
  It was
  shown by Seidenberg \cite{Se} that $D$ maps the normalization
  $k[\bar{S}]$ of $k[S]$ into itself.
  Then
  \[\sum s_jt^sf_j/t_j=D(t^s)\in k[\bar{S}]\sub k[t]
  \]
  implies $\prod_{j\not= i}t_j^{a_j-1}t^sf_i/t_i\in k[t]$.
  Hence $h_i$ divides
  $\prod_{j\not= i}t_j^{a_j-1}t^s/t_i$.
  But $t_i$ does
  not divide this product
  since $s_i=1$. Thus
  $h_i\in k$ and $f_i\in k[t]$. This means that $D$ restricts to a
  derivation of $k[t]$.
  \hfill $\Box$
\end{trivlist}
By Proposition 2 the Lie algebra $\TS=\Der k[S]$ of $k$-linear
derivations of the semigroup ring may be viewed as a subalgebra of
$\D=\Der k[t]$. Let us first describe the latter Lie algebra.
The derivations
$D_i=t_i\di$ span an Abelian subalgebra $H$. For a linear form $\l\in H^*$
let
\[\Dl=\{D\in\D,\;[h,D]=\l(h)\cdot D\;\mbox{for all}\; h\in H\}.
\]
Then $\D$ admits a root space decomposition
\[\D=\bigoplus_{\l\in H^*}\Dl.
\]
Given the basis
$D_1,\ldots,D_n$ of $H$ one may identify $H^*$ with $k^n$ by
identifying the form $\l$ with the vector
$(\l(D_1),\ldots,\l(D_n))$. Then the set of $\l\in H^*$
with $\Dl\not=0$ equals
\[\Nn\cup\{\l\in\Zn,\;\li=-1\;\mbox{for exactly one}\;i\;
  \mbox{and}\;\lj\ge0\;\mbox{for all}\;j\not= i\}.
\]
In fact, for $\l\in\Nn$ the root space $\Dl$ is spanned by all
$\Dlj=\tl t_j\dj$, $j=1,\ldots,n$. In particular, $\D_0=H$.
And if $\l\in\Zn$ with $\li=-1$
and $\lj\ge0$ for $j\not= i$ then
$\Dl$ is spanned by the single element $\Dli=\tl t_i\di$. All these
statements follow from the commutator relation
\[[D_i,\Dlj]=\li\cdot\Dlj.
\]
In order to describe the subalgebra $\TS$ we need some more notation.
Let
\begin{eqnarray*}
  \L_i & = & \{\l\in\Zn,\;  \l+s\in S\;\mbox{for all}\;
  s\in S \;\mbox{with}\; s_i\not= 0\},\;i=1,\ldots,n \\
  \L   & = & \L(S) = \bigcup \L_i \\
  \St  & = & \{\l\in\Nn,\; \l+(S\setminus\{0\})\sub S\}.
\end{eqnarray*}
\begin{trivlist}
  \item[] {\it Remarks.}
  (i) Let $n=1$. Then $k[S]$ is always Cohen-Macaulay, and the cardinality
  of $\L\setminus S$ equals the Cohen-Macaulay type of $k[S]$, see
  \cite{HK}. For $S=\N$ one has $\St=\N$ and $\L=\St\cup\{-1\}$.
  Otherwise $1\notin S$. Then our assumption that $\N\setminus S$ is
  finite implies $\L\sub\N$ and $\L=\St$.
  \item[]
  (ii) Let $n\ge2$. From $\l+\ai\in S$ for $\l\in\St$ and two indices
  $i$ one sees $\St\sub S'$. Hence $\St=S'$ in the Buchsbaum case and
  $\St=S$ in the Cohen-Macaulay case.
\end{trivlist}
\begin{trivlist}
  \item[] \bf Proposition 3. {\rm (i)} \it The Lie algebra $\TS$
  admits a root space decomposition
  \[\TS=\bigoplus_{\l\in H^*}\Tl.
  \]
  with $\Tl=\TS\cap\Dl$.
  \item[] {\rm(ii)} Suppose that $k[S]$ is Buchsbaum. Then
  the set of $\l\in H^*$ with $\Tl\not=0$ equals $\L(S)$. If $\l\in\St$
  then $\Tl$ is spanned by
  $D_{\lambda 1},\ldots,D_{\lambda n}$.
  And if $\l\in E_i=\L_i\setminus\St$ then
  $\Tl$ is spanned by the single element $\Dli$. In particular,
  $\L(S)=\St\cup\bigcup E_i$ is a disjoint union.
  \item[] \rm The elements of $\St$ (resp.\ $E_i$) will be called
  {\em ordinary} (resp.\ {\em i-exceptional}) {\em roots}.
  \item[] {\it Proof.} (i) For $D_\l=\sum_ib_{\l i}\Dli\in\Dl$ one has
  $D_\l t^s=\sum_ib_{\l i}s_i\cdot t^{\l+s}$. Hence $\sum_\l D_\l\in\TS$
  if and only if $\l+s\in S$ for all $s\in S$ and all occuring $\l$ with
  $\sum_i b_{\l i}s_i\not=0$ if and only if $D_\l\in\TS$
  for all occuring $\l$.
  \item[] (ii) Consider
  $\l\in\St$. Then
  $D_{\lambda 1},\ldots,D_{\lambda n}$
  are defined and contained in
  $\TS$. Next consider $\l\in \L_i$. From $\l+\ai\in S$ we see $\lj\ge0$
  for all $j\not= i$.
  Moreover, $\li\in\L(\pi_i(S))$ and Remark (i) above yields
  $\li\ge-1$.
  Hence $\Dli$ is defined and
  contained in $\TS$. Conversely, if $\Dli\in\TS$ then
  $\l\in \L_i$. The proof is completed by the following claim: If $\Tl$
  contains a linear combination of the
  $\Dli$ with at least two non-vanishing coefficients then $\l\in\St$.
  In fact,
  if $\sum_i b_i\Dli\in\TS$ with $b_1,b_2\not=0$ then $\l+\a^1$ and
  $\l+\a^2$ are contained in $S$. This gives $\l\in S'\sub\St$ as
  $k[S]$ is Buchsbaum.
  \hfill $\Box$
\end{trivlist}
\begin{trivlist}
  \item[] \it Examples. \rm
  (i) (\cite[Remark 1.3]{MT}) Let $S\sub\N^2$ be generated by
  (0,10),(3,7),(7,3),
  (8,2),(10,0) and let $\l=(9,11)$. Then $\l+(3,7)\notin S$ but $\l+s\in S$
  for the remaining generators $s$. Hence $\l\in S'\setminus\St$ and
  $k[S]$ is not Buchsbaum. Moreover, $\l\notin\L(S)$ but $\Tl\not=0$.
  In fact, $7D_{\l 1}-3D_{\l 2}\in\Tl$.
  \item[] (ii) Let $S\sub\N^2$ correspond to the affine cone over the
  $d$-uple embedding of $\P^1$ in $\P^d$, $d\ge2$, i.~e., $S$ is generated
  by $(0,d),(1,d-1),\ldots,(d-1,1),(d,0)$. Then $k[S]$ is normal and
  Cohen-Macaulay. The exceptional roots are $(-1,1)+m(0,d)$ and
  $(1,-1)+m(d,0)$ with $m\in\N$.
  \item[] (iii) Let $S\sub\N^2$ correspond to the product of a cusp with
  a line, i.~e., $S$ is generated by $(2,0)$, $(3,0)$ and $(0,1)$.
  Then $k[S]$
  is Cohen-Macaulay. The 1-exceptional roots are $(1,0)+m(0,1)$ with
  $m\in\N$. The 2-exceptional roots are $(0,-1)+m(2,0)$ and $(3,-1)+m(2,0)$
  with $m\in\N$.
\end{trivlist}
Examples (ii) and (iii) illustrate the second part of the next result.
\begin{trivlist}
  \item[] \bf Proposition 4. \it {\rm (i)}
  $\St$ is a finitely generated subsemigroup of
  $\Nn$.
  \item[] {\rm (ii)} Suppose that $k[S]$ is Buchsbaum and $n\ge2$.
  For fixed $i$ let
  $A_i$ be the semigroup generated by all $\a^j$ with $j\not= i$. Then
  the set $E_i$ of $i$-exceptional roots is a finitely generated
  $A_i$-module.
  \item[] Proof. \rm
  (i) Clearly $\St$ is a subsemigroup of $\Nn$. Let $A$ be the semigroup
  generated by $\a^1,\ldots,\a^n$. We show more generally that every
  subsemigroup $T\sub\Nn$ containing $A$ is finitely generated. Let $a_i$
  be the nonzero entry of $\ai$. For $\b\in\Nn$ with $\b_i<a_i$ for all $i$
  let $T_\b=(\b+A)\cap T$. By Dickson's Lemma
  each $T_\b$ is a finitely generated $A$-module (or empty).
  Since $T=\bigcup T_\b$ is a finite union, $T$ is
  finitely generated as an $A$-module and hence as a semigroup.
  \item[] (ii) We may assume $i=1$. If $\l\in E_1=\L_1\setminus\St$
  then clearly
  $\l+\a^2\in\L_1$. Moreover, $\l+\a^1\in S$ so that $\l\in S'_i$ for
  $i\ge2$. If $\l+\a^2\in\St$ then $\l+2\a^2\in S$, hence $\l\in S'_1$
  and $\l\in S'=\St$, contradiction. Thus $\l+\a^2\in E_1$. This proves
  that $E_1$ is an $A_1$-module. It remains to show that it is finitely
  generated. For $\gamma\in\N\times\{0\}\sub\Nn$ and
  $\b\in\{0\}\times\N^{n-1}\sub
  \Nn$ with $\b_i<a_i$ for all $i$ let
  $E_{\gamma\b}=(\gamma+\b+A_1)\cap E_1$. As above this is a finitely
  generated $A_1$-module (or empty). If $E_{\gamma\b}\not=\emptyset$ and
  $\gamma'=\gamma+m\a^1$ for some $m\in\N$, $m\not=0$ then
  $E_{\gamma'\b}=\emptyset$. Otherwise, there is $\l\in A_1$ with
  $\gamma+\b+\l,\gamma'+\b+\l\in E_1$, contradicting
  $\gamma'+\b+\l=\gamma+\b+\l+m\a^1\in S\sub\St$.
  Since there are only finitely many congruence classes
  of $\N$ modulo $\a^1$ the Proposition is proven.
  \hfill $\Box$
\end{trivlist}
\section{Reconstruction of the semigroup}
Before we explain how to reconstruct the semigroup $S$ from its
Lie algebra $\TS$ we make a remark concerning the reconstruction of $S$
from its semigroup ring $k[S]$ discussed by Gubeladze \cite{Gu}.
Consider the augmentation $k[S]\to k$ defined by $t^s\mapsto 0$ for
all $s\in S\setminus\{0\}$. Gubeladze \cite[Theorem 2.1]{Gu} proved
that affine semigroups $S_1$ and $S_2$ are isomorphic if $k[S_1]$
and $k[S_2]$ are isomorphic as augmented algebras. Moreover
\cite[Lemma 2.8]{Gu}, if $k[S_1]$ and $k[S_2]$ are normal and isomorphic
just as algebras then they are isomorphic as augmented algebras. We shall
extend this result (for simplicial semigroups) to the Buchsbaum case.
\\[1ex]
Let us say that $S$ corresponds to a {\em product along a line} if, after
permutation of coordinates, $S=\N\oplus M$ for some semigroup
$M\sub\N^{n-1}$. We shall see that this property only depends on the
algebra $k[S]$ and even on the Lie algebra $\TS$.
Let $L=[\TS,\TS]$ be the derived algebra.
\begin{trivlist}
  \item[] \bf Proposition 5. \it Suppose that $k[S]$ is Buchsbaum.
  Then the following are equivalent:
  \begin{itemize}
    \item[{\rm (a)}] The semigroup $S$ corresponds to a product along a
    line.
    \item[{\rm (b)}] There is $\l\in\L(S)$ with $|\l|<0$.
    \item[{\rm (c)}] $L=\TS$.
  \end{itemize}
  \item[] Proof. \rm (a) $\Leftrightarrow$ (b) If $(-1,0,\ldots,0)$ is a
  root then $(1,0,\ldots,0)\in S$ and $S=\N\oplus M$
  with $M=S\cap\ker\:\pi_1$. The
  converse is clear.
  \item[] (b) $\Rightarrow$ (c) Here and later we use the commutator
  relation
  \[[\Dli,\Dmj]=\mi\Dlmj-\lj\Dlmi.
  \]
  It shows
  $\bigoplus_{\l\not=0}\Tl\sub L$.
  Let $\l=(-1,0,\ldots,0)\in\L$ so that
  $\m=(1,0,\ldots,0)\in S\sub\St$. Then $L$ contains
  $2D_1=[D_{\l 1},D_{\m 1}]$ and $D_j=[D_{\l 1},D_{\m j}]$ for $j\ge2$.
  Thus $\T_0=H\sub L$.
  \item[] (c) $\Rightarrow$ (b) Assume that $|\l|\ge0$ for all roots $\l$.
  Then $\eta^1+\eta^2=0$ for roots $\eta^1,\eta^2\not=0$
  is possible only if (after permutation
  of coordinates) $\eta^1=(-1,1,0,\ldots,0)$, $\eta^2=(1,-1,0,\ldots,0)$.
  In this case $[D_{\eta^1,1},D_{\eta^2,2}]=D_2-D_1$. Since $\T_0$ is
  Abelian we obtain
  \[L\sub\bigoplus_{\l\not=0}\Tl\oplus<D_n-D_1,\ldots,D_2-D_1>
  \]
  and $\T_0\not\sub L$.
  \hfill $\Box$
\end{trivlist}
\begin{trivlist}
  \item[] \bf Proposition 6. \it
  Suppose that $k[S_1]$ and $k[S_2]$ are Buchsbaum.
  \item[] {\rm (i)} If $k[S_1]$ and $k[S_2]$ are isomorphic as algebras
  then they are isomorphic as augmented algebras.
  \item[] {\rm (ii)} If $S_1$ and $S_2$ do not correspond to products
  along a line then every algebra isomorphism $\phi:k[S_1]\to k[S_2]$
  is augmented.
  \item[] Proof. \rm
  Let $I\sub k[S_2]$ be a proper differential ideal, i.~e.,
  $D(I)\sub I$ for every $D\in\T(S_2)$. We claim that $I$ is generated
  by some monomials $t^s$, $s\in S_2$. In particular, $I$ is contained
  in the augmentation ideal generated by all $t^s$,
  $s\in S_2\setminus\{0\}$. Given $f=\sum b_st^s\in I$ fix any $s$
  with $b_s\not=0$. Take any of the remaining $\l\in S_2$ with
  $b_\l\not=0$ and choose $j$ with $\l_j\not= s_j$. Then
  $\sum_\m(\lj-\m_j)b_\m t^\m=\lj f-D_j(f)\in I$
  contains less monomials than $f$ but still the monomial $t^s$.
  Repeated application yields $t^s\in I$, proving the claim.
  \\[1ex]
  Now assume $S_1=\N^m\oplus M$ for some $M\sub \N^{n-m}$ which does
  not correspond to a product along a line. Let $J$ be the ideal
  of $k[S_1]$ generated by all $t^\m$, $\m\in M\setminus\{0\}$.
  We claim that $J$ is differential. Consider any $\l\in\L_i$,
  $i=1,\ldots,n$. In order to show $\Dli(t^\m)=\mi t^{\l+\m}\in J$ we
  may assume $\mi\not=0$. Then $\l+\m\in S_1$. From $|\m|\ge2$ we
  conclude $\l+\m=\n+\m'$ with $\n\in\N^m$ and $\m'\in M\setminus\{0\}$.
  Hence $t^{\l+\m}=t^{\n+\m'}\in J$.
  \\[1ex]
  Let $\phi:k[S_1]\to k[S_2]$ be an algebra isomorphism. It induces
  a Lie algebra isomorphism $\phi^\sharp:\T(S_1)\to\T(S_2)$ by
  $D\mapsto\phi\circ D\circ\phi^{-1}$. Since $J$ is differential
  its image in $k[S_2]$ is differential and hence contained in the
  augmentation ideal of $k[S_2]$. We have $k[S_1]=k[M][t_1,\ldots,t_m]$.
  For $i=1,\ldots,m$ let $c_i$ be the constant term of $\phi(t_i)$.
  Define the $k[M]$-automorphism $\psi$ of $k[S_1]$ by
  $\psi(t_i)=t_i-c_i$, $i=1,\ldots,m$. Then the $\phi\circ\psi(t_i)$
  have no constant term. Since the augmentation ideal of $k[S_1]$
  is generated by $t_1,\ldots,t_m$ and $J$ this means that
  $\phi\circ\psi$ is augmented. Assertion (ii) now also is clear
  because in that case $J$ equals the augmentation ideal.
  \hfill $\Box$
\end{trivlist}
\begin{trivlist}
  \item[] \bf Theorem 1. \it Let $S_1,S_2$ be simplical affine
  semigroups such that $k[S_1],k[S_2]$ are Buchsbaum. Suppose that the
  Lie algebras $\T(S_1),\T(S_2)$ are isomorphic. Then $S_1,S_2$ have the
  same rank and the semigroups $\tilde{S_1},\tilde{S_2}$ are isomorphic.
  \item[] Proof. \rm
  If $\T(S_1)$ equals its derived algebra then $S_1$ and $S_2$
  correspond to products along a line. By a result of Skryabin
  \cite[Theorem 2]{Sk} the semigroup rings $k[S_1],k[S_2]$ are isomorphic.
  Then \cite[Theorem 2.1]{Gu} and Proposition 6 imply that the semigroups
  $S_1,S_2$ themselves
  are isomorphic. Now suppose that the derived algebra is
  strictly smaller than $\T(S_1)$. Then $|\l|\ge 0$ for all
  $\l\in\L(S_1)$. As $[\Tl,\Tm]\sub\Tlm$ for all roots $\l,\m$ the
  subspaces $I_d=\bigoplus_{|\l|\ge d}\Tl$ are ideals of $\T(S_1)$
  with finite dimensional quotients $\T(S_1)/I_d$ and
  $\bigcap_{d\in\N}I_d=0$. Given an isomorphism $\T(S_1)\simeq\T(S_2)$
  we obtain an Abelian subalgebra $H_2$ of $\T(S_1)$ and another
  root space decomposition $\T(S_1)=\bigoplus_{\m\in H_2^*}\Tm'$.
  Every finite dimensional subspace of $\T(S_1)$ is mapped isomorphically
  onto its image in $\T(S_1)/I_d$ if $d$ is sufficiently large. Thus,
  for $d\gg 0$, $H_2$ embeds into $Q=\T(S_1)/I_d$.
  For $\m\in H_2^*$ consider the root spaces
  \[Q_{\m}'=\{D\in Q,\; [h,D]=\m(h)\cdot D\;\mbox{for all}\;h\in H_2\}.
  \]
  Their sum is direct. Since each $\Tm'$ is mapped into
  $Q_{\m}'$ and the images of the $\Tm'$ span $Q$ we see
  $Q=\bigoplus_{\m\in H_2^*}Q_{\m}'$ and that each $\Tm'$ is
  mapped onto $Q_{\m}'$. In particular, $Q_0'=H_2$.
  It follows that $H_2$ equals its normalizer in $Q$
  and hence is a Cartan subalgebra of $Q$.
  Using Proposition 3, Remark (i) preceding it,
  and Proposition 4 we may assume that
  the subsemigroup of $H_2^*$ generated by all $\m$ with
  $\dim\:Q_{\m}'=\dim\:H_2=\rk S_2$ equals $\tilde{S_2}$.
  Analogous statements hold true for $H_1$ and $d\gg 0$. Since $Q$ is
  finite dimensional there is an automorphism of $Q$ mapping the
  Cartan subalgebra $H_1$ onto the second Cartan subalgebra $H_2$,
  \cite[section 16]{Hu}. Its dual induces an isomorphism between
  the semigroups $\tilde{S_1}$ and $\tilde{S_2}$.
  \hfill $\Box$
\end{trivlist}
Using Remark (ii) preceding Proposition 3 we conclude
\begin{trivlist}
  \item[] \bf Corollary 1. \it
  Simplicial affine semigroups $S$ of rank $\ge 2$ with $k[S]$
  Cohen-Macaulay are uniquely determined by their Lie algebra $\TS$.
\end{trivlist}
Look again at Gubeladze's Theorem that $S$ is uniquely determined
by the augmented algebra
$k[S]$. In the above proof we applied this only in case $S$ does
correspond to a product along a line. Therefore,
using the Lie algebra $\TS$ as an intermediate step,
we have reproven Gubeladze's Theorem
in the special case that $S$ is simplicial,
does not correspond to a product along a line, and
$k[S]$ is Cohen-Macaulay of dimension $\ge 2$. But $\TS$ cannot
distinguish between semigroups with the same Cohen-Macaulayfication:
\begin{trivlist}
  \item[] \it Examples. \rm
  (i) Fix $d,l\in\N$, both $\ge2$. Let $S$ consist of all $s\in\N^2$ with
  $|s|=md$, $m\ge l$. Then $k[S]$ is Buchsbaum and the
  Cohen-Macaulayfication $S'$ is generated by
  $(0,d),(1,d-1),\ldots,(d-1,1),(d,0)$. Both $S$ and $S'$ have the same
  exceptional roots, see Example (ii) after Proposition 3. Hence
  $\TS=\T(S')$, independently of $l$.
  \item[] (ii) Let $S_1$ (resp.\ $S_2$) be generated by all $\l\in\N^2$
  with $|\l|=6$ except $\l=(3,3)$ (resp.\ $\l=(2,4)$).
  They have a Buchsbaum semigroup ring and the same
  Cohen-Macaulayfication generated by all $\l\in\N^2$ with $|\l|=6$.
  In both cases the exceptional roots are $(-1,7)+m(0,6)$ and
  $(7,-1)+m(6,0)$ with $m\in\N$. Hence $\T(S_1)=\T(S_2)$. But $S_1,S_2$
  are not isomorphic. In fact, any isomorphism would map the set of extremal
  elements $\{(6,0),(0,6)\}$ onto itself, hence
  $(6,6)$ onto $(6,6)$. This contradicts $(6,6)=2(3,3)$ in $S_2$ but
  $(6,6)\not=2s$ for all $s\in S_1$. Observe that both semigroups
  correspond to affine cones over smooth projective curves in $\P^5$.
\end{trivlist}
In the rank 1 case the situation is different. Although the semigroup
ring always is Cohen-Macaulay the semigroup is, in general, not determined
by the Lie algebra:
\begin{trivlist}
  \item[] \it Examples. \rm
  (i) The numerical semigroups generated by 2 and 3 (resp.\ 3, 4 and 5)
  have the same $\St=\N$, hence the same Lie algebra. Observe that the
  semigroup ring is Gorenstein in the first case whereas it has
  Cohen-Macaulay type 2 in the second, see Remark (i) preceding
  Proposition 3.
  \item[] (ii) The numerical semigroups generated by 3, 7 and 8
  (resp.\ 4, 5 and 7) have the same $\St$ generated by 3, 4 and 5,
  hence the same Lie algebra. Observe that the Cohen-Macaulay type is 2
  in both cases.
\end{trivlist}
\begin{trivlist}
  \item[] \bf Corollary 2. \it Numerical semigroups $S$ with $k[S]$
  Gorenstein
  are uniquely determined by $\TS$ and even by the finite dimensional
  Lie algebra $\TS/[L,L]$.
  \item[] Proof. \rm
  If $L=\TS$ then $S=\N$. So suppose $L\not=\TS$.
  Then $\St$ is the set of roots and $L=\bigoplus_{\l\not=0}\Tl$.
  This implies
  $\Tl\cap[L,L]=0$
  for $\l$ in the minimal generator
  system of $\St$ and $\Tl\sub[L,L]$ for every $\l$ which can be decomposed
  as $\l=\m+\n$ with two different $\m,\n\in\St$.
  We see that $\TS/[L,L]$ is finite dimensional and that we can use
  the intrinsically defined ideal $[L,L]$ instead of $I_d$ in the proof
  of Theorem 1.
  It remains to show that $S$ is uniquely determined
  by $\St$ in the Gorenstein case. By \cite[Satz 1.9, Proposition 2.21]{HK}
  we know
  $\St=S\cup\{c-1\}$ with the conductor $c$ of $S$.
  Consider first the case $\St=\N$. Then $S$ must be the semigroup
  $\N\setminus\{1\}$, generated by 2 and 3. Now let $\St\not=\N$.
  Let $a$ be the
  smallest element of $S$ different from $0$. As $S$ is a symmetric
  semigroup we see $c-2,\ldots,c-a\in S$ but $c-a-1\notin S$. Thus $\St$
  has conductor $c-a$.
  Then $c-a\in S\setminus\{0\}$ implies
  $c-1> c-a\ge a$. Hence $a$ is the smallest element of $\St$ different
  from $0$. Therefore, $S=\St\setminus\{c-1\}$ is determined via $c-a$
  and $a$ by $\St$.
  \hspace*{\fill} $\Box$
\end{trivlist}
\section{Automorphisms of the Lie algebra}
Every automorphism $\phi$ of $k[S]$ induces a Lie algebra automorphism
\[\phi^{\sharp}:\TS\to\TS:D\mapsto\phi\circ D\circ\phi^{-1}.
\]
The purpose of this section is to show
\begin{trivlist}
  \item[] \bf Theorem 2. \it
  Let $S$ be a simplicial affine semigroup such that $k[S]$ is
  Cohen-Macaulay. For every automorphism $\Phi$ of $\TS$ there is a unique
  automorphism $\phi$ of $k[S]$ such that $\Phi=\phi^{\sharp}$.
  \item[] Proof. \rm
  If $\Phi=\phi^{\sharp}$ then
  $\Phi(f\cdot\Phi^{-1}(D))=\phi(f)\cdot D$ for all $f\in k[S]$ and
  $D\in\TS$. This shows uniqueness. Now take an arbitrary automorphism
  $\Phi$ of $\TS$. If $S$ corresponds to a product along a line the
  assertion follows from \cite[Theorem 2]{Sk}. By Proposition 5
  we may assume $|\l|\ge 0$ for all
  $\l\in\L$. The Lie algebra $\TS$ is graded by
  $\T^d=\bigoplus_{|\l|=d}\Tl$. Note that $\T^0$ consists of the linear
  vector fields in $\TS$, i.~e., those $\sum f_i\di\in\TS$ where the $f_i$
  are linear forms in the variables $t_i$.
  The homogeneous component of smallest degree
  of $D\in\TS$, $D\not=0$,
  will be called the leading
  form of $D$.
  We claim that every $h'=\Phi(h)\in\Phi(H)$, $h'\not=0$,
  has leading form of degree zero.
  In fact, choose $\l\in\L$ with $\l(h)\not=0$ and $Y\in\Phi(\Tl)$,
  $Y\not=0$. Comparison of leading forms in
  $[h',Y]=\l(h)\cdot Y$ yields the claim.
  Hence the leading forms of the vector fields $Y_i=\Phi(D_i)$,
  $i=1,\ldots,n$, are linear vector fields and linearly independent.
  We can thus find a point $p$ in affine space $\A^n$ such that the
  tangent vectors $Y_1(p),\ldots,Y_n(p)$ are linearly independent.
  Now consider the polynomial ring $k[t]$ as a subring of the ring
  $\F=k[[t-p]]$ of formal power series centered at $p$ and $\Der k[t]$
  as a subalgebra of the Lie algebra $\Der \F$. By Proposition 7
  below there are formal coordinates
  $s_1,\ldots,s_n$ at $p$, i.~e., elements of $\F$
  vanishing at $p$ with $k[[s]]=\F$,
  such that $Y_i=\d_{s_i}$ in $\Der \F$.
  Let $x_i=\exp\; s_i$ for $i=1,\ldots,n$. If $\l\in\Zn$ then
  \[Y_i(x^\l)=\li\cdot x^\l \quad\mbox{for}\quad i=1,\ldots,n
  \]
  and, up to multiplication with a constant, $x^\l$ is the unique
  element of $\F$ with this property. This implies that for $\l\in\Zn$
  the root space
  \[\Tl'=\{D\in\Der\F,\;[Y_i,D]=\li\cdot D\;\mbox{for all}\;i\}
  \]
  is spanned by the $\Yli=x^\l Y_i$, $i=1,\ldots,n$. We conclude that
  $\Phi(\Tl)=\Tl'$ for ordinary roots $\l\in\St$ and
  $\Phi(\Tl)\sub\Tl'$ for exceptional roots $\l\in\bigcup E_i$. Next we
  claim
  \[\Phi(\Dli)=b_{\l i}\Yli \quad\mbox{for all $\l$ and $i$}
  \]
  with suitable constants $b_{\l i}\not=0$. To prove this, note that
  $[\Dli,\Dmj]=-\lj \Dlmi$ if $\mi=0$ and thus
  $Y=\Phi(\Dli)$ has the following property:
  For all $\m\in\St$ with $\mi=0$ the image of
  $\ad Y: \Tm'\to\Tlm'$ has dimension $\le 1$. Hence it is enough to
  show that, up to multiplication with a constant, $\Yli$ is the unique
  element of $\Tl'$ with this property. In fact, for
  $Y=\sum_k c_kY_{\l k}$ the matrix of coefficients of $([Y,Y_{\m j}])_j$
  with respect to the basis $(Y_{\l+\m,k})_k$ has determinant equal to
  the value at $\sum c_k\m_k$ of the characteristic polynomial of
  the matrix $(\lj c_k)_{j,k}$. This value does not vanish
  for a suitable
  choice of $\m\in\St$ with $\mi=0$ if $c_k\not=0$ for some $k\not= i$,
  and the claim is proven.
  \\[1ex]
  For fixed $\l$ choose $\m\in\St$ with $\m_1\not=\l_1$ and
  $\mi\not=0$ for $i\not=1$. Then the usual commutator relation
  implies
  $b_{\l i}b_{\m 1}=b_{\l+\m,1}$ for all $i$.
  Hence the $b_{\l i}$ are independent of $i$, say
  $b_{\l i}=b_\l$. We have
  \[\Phi(\Dli)=b_\l\Yli \quad\mbox{for all $\l$ and $i$}.
  \]
  Denote by $\G$ the subgroup of $\Zn$ generated by $\L$.
  As $b_\l b_\m=b_{\l+\m}$ for all $\l,\m\in\L$
  with $\l+\m\in\L$ the map
  $\l\mapsto b_\l$ can be extended
  to a homomorphism $\G\to k^*$. The group $\G$ is
  free of rank $n$, say generated by $\gamma_1,\ldots,\gamma_n$.
  There is a rational matrix $Q=(q_{ij})_{i,j}$ such that
  $l=Q\cdot\l$ if $\l=\sum l_i\gamma_i\in\G\sub\Zn$ and
  $l=(l_1,\ldots,l_n)\in\Zn$. Write $q_{ij}=r_{ij}/s$ with integers
  $r_{ij},s$, choose $\zeta_i\in k$ such that $b_{\gamma_i}=\zeta_i^s$,
  and let $c_j=\prod_i \zeta_i^{r_{ij}}$. Then
  $b_\l=c^\l$ for all $\l\in\G$.
  Thus, if we replace the $x_j$ by $c_jx_j$ we obtain
  for the new $\Yli=x^\l Y_i$ the equations
  \[\Phi(\Dli)=\Yli\quad\mbox{for all $\l$ and $i$}.
  \]
  We have seen above
  that $\TS$ is spanned by all $x^\l Y_i$ with $\l\in\L_i$ and
  $i=1,\ldots,n$.
  Using $|\l|\ge0$
  for all $\l$ it is easy to show that each $\L_i$ is an $\St$-module.
  Hence $\TS$ is a
  module over the subalgebra of $\F$ generated by all $x^s$, $s\in\St$.
  Fix $s\in\St$. From
  $x^st_i\d_i=x^sD_i\in\TS\sub\Der k[t]$ we conclude that the element
  $x^st_i$ of $\F$ actually is contained in the subalgebra $k[t]$.
  Since the same is true for $x^{2s}t_i$ we obtain $x^s\in k[t]$.
  Even more: $x^sD_i\in\TS$ for $i=1,\ldots,n$ shows $x^s\in k[\St]$.
  The $x_i$ are algebraically independent. Therefore, an algebraic
  relation between finitely many $t^{s_1},\ldots,t^{s_m}$ holds if and
  only if the same relation holds between
  $x^{s_1},\ldots,x^{s_m}$. This means that we can define an injective
  homomorphism $\phi:k[\St]\to k[\St]$ by $\phi(t^s)=x^s$.
  The equations $\Phi(t^\l D_i)=x^\l Y_i$ established above translate
  into $\Phi(D)\circ\phi=\phi\circ D$ for all $D\in\TS$.
  Using $\Phi^{-1}$ instead of $\Phi$ we get an injective endomorphism
  $\psi$ of $k[\St]$ with $D\circ\psi=\psi\circ\Phi(D)$ for all $D$.
  Then $D_i\circ\psi\circ\phi=\psi\circ\phi\circ D_i$ for all $i$. Using
  this information one shows that $\psi\circ\phi$ maps each
  one-dimensional subspace of $k[\St]$ spanned by some $t^s$
  into itself. Hence injectivity of $\phi$ and $\psi$ implies surjectivity
  of both. In case $n\ge2$ we are done because then $S=\St$
  by our hypothesis on $k[S]$ and $\phi$ is an
  automorphism of $k[S]=k[\St]$ with $\Phi(D)=\phi\circ D\circ\phi^{-1}$
  for all $D\in \TS$, i.~e., $\Phi=\phi^\sharp$.
  \\[1ex]
  Finally, consider the case $n=1$. Then $x$ is a single element of $\F$
  with $x^s\in k[t]$ for all $s\in\St$. Since $\St$ is a numerical
  semigroup $x$ must be contained in $k(t)$ and, being integral over
  $k[t]$, even in $k[t]$. As $\phi:t^s\mapsto x^s$ defines an
  automorphism of $k[\St]$ the polynomial $x$ has degree 1, say
  $x=a+bt$. We had $Y(x)=x$ for $Y=\Phi(t\d_t)\in\TS$. This implies
  $Y=(a/b+t)\d_t$ and $a=0$ because $S\not=\N$, i.~e., $-1$ is not a root.
  Therefore, $\phi$ restricts to an automorphism of $k[S]$ with
  $\Phi=\phi^\sharp$.
  \hspace*{\fill} $\Box$
\end{trivlist}
It remains to show
\begin{trivlist}
  \item[] \bf Proposition 7. \it
  Let $\F=k[[t_1,\ldots,t_n]]$. Suppose that $Y_1,\ldots,Y_n\in\Der \F$
  satisfy $[Y_i,Y_j]=0$ for all $i,j$ and that $Y_1(0),\ldots,Y_n(0)$
  are linearly independent. Then there are formal coordinates
  $s_1,\ldots,s_n$ such that $Y_i=\d_{s_i}$ for all $i$.
  \item[] Proof. \rm
  Write $Y_i=\sum_jf^i_j\d_{t_j}$. By hypothesis the matrix
  $F=(f^i_j)_{i,j}$ is invertible over $\F$, say with inverse
  $G=(g^j_k)_{j,k}$. Application of $Y_m$ to
  \[\sum_jf^i_jg^j_k=\delta^i_k
  \]
  yields
  \renewcommand{\theequation}{\fnsymbol{equation}}
  \begin{equation}
    \sum_{l,j}f^m_lf^i_j(\d_{t_l}g^j_k)
    =-\sum_{l,j}f^m_l(\d_{t_l}f^i_j)g^j_k
  \end{equation}
  The hypothesis $[Y_m,Y_i]=0$ means
  \[\sum_lf^m_l(\d_{t_l}f^i_j)=\sum_lf^i_l(\d_{t_l}f^m_j)
  \]
  for all $j$. Hence we may interchange $i$ and $m$ in the right hand
  side and, therefore, in the left hand side
  of (*). After renaming the summation indices we obtain
  \[\sum_{l,j}f^m_lf^i_j(\d_{t_l}g^j_k-\d_{t_j}g^l_k)=0.
  \]
  Invertibility of $F$ implies
  \[\d_{t_l}g^j_k=\d_{t_j}g^l_k
  \]
  for all $l$, $j$ and $k$. This condition is equivalent
  (over a field of characteristic 0) to the existence of
  $s_1,\ldots,s_n\in\F$ vanishing at 0 with
  \[g^j_k=\d_{t_j}s_k
  \]
  for all $j$ and $k$. These $s_k$ form  a system of coordinates because
  $G$ is invertible. And clearly $Y_is_k=\delta^i_k$ for all $i,k$.
  \hfill $\Box$
  \renewcommand{\theequation}{\arabic{equation}}
  \setcounter{equation}{0}
\end{trivlist}
\begin{trivlist}
  \item[] \it Remark. \rm
  Proposition 7 is the special case $r=n$ of a more general statement
  involving an arbitrary number $r\le n$ of vector fields. The latter
  usually is stated for differentiable or analytic vector fields
  over the fields of real or complex numbers and appears in the
  literature in connection with Frobenius' Theorem. It is surely known
  to hold for formal power series vector fields over arbitrary fields
  of characteristic 0. But lacking an explicit reference we have chosen
  to provide the very simple proof above.
\end{trivlist}
\section{Derivations of the Lie algebra}
In this section we show
\begin{trivlist}
  \item[] \bf Theorem 3. \it Let $S\sub\Nn$ be a simplicial affine
  semigroup such that $k[S]$ is Buchsbaum.
  Then every derivation $\DE$ of $\TS$ is inner: $\DE=\ad D$ for some
  $D\in\TS$.
  \item[] Proof. \rm
  The cochain complex of the Lie algebra $\TS$ with coefficients in the
  adjoint representation has a
  $\Zn$-grading given by the root space decomposition.
  By \cite[Theorem 1.5.2b]{F} it is acyclic in degrees different from zero.
  Hence we may assume that
  the given $\DE$ has degree $0$, i.~e.\ $\DE(\Tl)\sub\Tl$ for all
  $\l$.
  For each root $\l$ denote by $M(\l)$ the set of $i$ such that $\Dli\in\TS$.
  Thus $M(\l)=\{1,\ldots,n\}$ for ordinary roots
  and $M(\l)=\{i\}$ for $i$-exceptional roots.
  We have
  \begin{equation}
    \DE(\Dli)=\sum_{m\in M(\l)}b_{\l im}D_{\l m}
    \quad\mbox{for}\quad i\in M(\l)
  \end{equation}
  with suitable constants $b_{\l im}\in k$. The brackets of the generators
  are given by
  \begin{equation}
    [\Dli,\Dmj]=\mi \Dlmj-\lj \Dlmi
  \end{equation}
  Inserting (1) and (2) into the cocycle condition
  \[\DE([\Dli,\Dmj])=[\DE(\Dli),\Dmj]+[\Dli,\DE(\Dmj)]
  \]
  gives
  \begin{eqnarray*}
    \lefteqn{
    \sum_m(\mi\cdot b_{\l+\m,j,m}-\lj\cdot b_{\l+\m,i,m})D_{\l+\m,m}}
    \hspace{2cm}
    \\ & = &
    \sum_m(\mi\cdot b_{\m jm}-\lj\cdot b_{\l im}) D_{\l+\m,m}
    \\ &   &
    +(\sum_m\m_m\cdot b_{\l im})\Dlmj
    -(\sum_m\lm\cdot b_{\m jm})\Dlmi.
  \end{eqnarray*}
  By comparing the coefficients one obtains
  \begin{equation}
    \mi\cdot b_{\l+\m,j,m}-\lj\cdot b_{\l+\m,i,m}=
    \mi\cdot b_{\m jm}-\lj\cdot b_{\l im}
    \quad \mbox{for} \quad m\not=i,j
  \end{equation}
  \begin{equation}
    \mi\cdot b_{\l+\m,j,j}-\lj\cdot b_{\l+\m,i,j}=
    \mi\cdot b_{\m jj}-\lj\cdot b_{\l ij}
    +\sum_m\m_m\cdot b_{\l im}
    \quad \mbox{for} \quad j\not=i
  \end{equation}
  \begin{equation}
    (\mi-\li)b_{\l+\m,i,i}=\mi\cdot b_{\m ii}-\li\cdot b_{\l ii}
    +\sum_m\m_m\cdot b_{\l im}-\sum_m\lm\cdot b_{\m im}
  \end{equation}
  Equation (4) with $\l=\m=\a^j$ yields
  \begin{equation}
    b_{2\a^j,i,j}=0\quad \mbox{for} \quad i\not=j
  \end{equation}
  Let us show that $b_{\l ij}=0$ for all $\l\in\St$
  and all $i,j\in M(\l)$ with $i\not= j$.
  Set $\m=2\a^j$. In case $\li=0$ the claim follows from (5) and (6).
  If $\li\not=0$ use (3) with $j=i$ and $m$ replaced by $j$ to show
  $b_{\l+\m,i,j}=b_{\l ij}$. Then (4) gives the claim.
  \\[1ex]
  Now we have
  \[\DE(\Dli)=b_{\l i}\Dli\quad\mbox{for}\quad i\in M(\l)
  \]
  with suitable $b_{\l i}\in k$. Equations (4) and (5) reduce to
  \begin{equation}
    \mi\cdot b_{\l+\m,j}=\mi\cdot b_{\m j}+\mi\cdot b_{\l i}
    \quad\mbox{for}\quad j\not=i
  \end{equation}
  \begin{equation}
    (\m_j-\lj)b_{\l+\m,j}=(\m_j-\lj)(b_{\l j}+b_{\m j})
  \end{equation}
  For fixed $\l\in\St$ the coefficients
  $b_{\l i}$ are independent of $i\in M(\l)$. In fact, for $j\not=i$
  apply (7) and (8) where
  $\m$ is any element of $\St$ with $\mi\not=0$ and $\m_j\not=\lj$.
  Thus we may write $b_\l$ instead of $b_{\l i}$.
  \\[1ex]
  Consider first the case $n\ge2$. Then (7) implies
  $b_{\l+\m}=b_\l+b_\m$ for $\l,\m\in\St$. Let $c_i=b_{\ai}/a_i$ where
  $a_i$ denotes the nonzero entry of $\ai$.
  Using the fact that $\St$ is torsion modulo the semigroup generated
  by the $\ai$ one shows
  $b_\l=\sum_ic_i\l_i$ for all $\l\in\St$. The same is seen to hold
  for $\l\in\L_i$ by applying (7) with some $\m\in S$, $\mi\not=0$.
  We have proven
  \[[\sum_ic_i D_i,\Dlj]=\sum_ic_i\li\Dlj=b_\l\Dlj=\DE(\Dlj)
  \]
  for all $\l\in\L$ and $j\in M(\l)$. This means $\DE=\ad D$ for
  $D=\sum_ic_iD_i$.
  \\[1ex]
  In the case $n=1$ only equation (8) is available. Then
  $b_{5\l}=b_{3\l}+b_{2\l}=2b_{2\l}+b_\l$ and
  $b_{5\l}=b_{4\l}+b_\l=b_{3\l}+2b_\l=b_{2\l}+3b_\l$, hence
  $b_{2\l}=2b_\l$ and then $b_{m\l}=m b_\l$ for all $m\in\N,\l\in\St$ with
  $m,\l>0$. This shows that the ratio $b_\l/\l$ is independent of $\l$,
  say $b_\l/\l=c$. Hence $b_\l=c\l$ for all positive roots. Since the same
  clearly holds for $\l=0$ (and $\l=-1$ in the special case $S=\N$) we
  have again shown that $\DE$ is inner.
  \hfill $\Box$
\end{trivlist}
\begin{trivlist}
  \item[] \it Remark. \rm In the special case $S=\Nn$ Theorem 2 was
  proven by Heinze \cite[Kap.~II, Satz 2.8]{He}. More generally,
  for semigroups corresponding to a product along a line it follows
  from work of Skryabin \cite[Theorem 3]{Sk}.
\end{trivlist}

\vfill
Departamento de Algebra, Geometria y Topologia, Universidad de Valladolid \\
E 47005 Valladolid, Spain \\
campillo@cpd.uva.es \\[2ex]
Instytut Matematyki, Uniwersytet Warszawski \\
PL 02-097 Warszawa, Poland \\
jagrab@mimuw.edu.pl \\[2ex]
Fachbereich Mathematik, Universit\"at Mainz \\
D 55099 Mainz, Germany \\
mueller@mat.mathematik.uni-mainz.de
\end{document}